\documentclass[smallextended,referee,envcountsect]{svjour3} 
\smartqed 
\usepackage{graphicx}
\journalname{JOTA}
\usepackage{mathtools}
\usepackage[normalem]{ulem} 

\usepackage{amssymb}
\usepackage{mathrsfs}
\usepackage{color}

\newcommand{\overbar}[1]{\mkern 1.5mu\overline{\mkern-1.5mu#1\mkern-1.5mu}\mkern 1.5mu}

\newcommand{\cut}[1]{{}}

\newcommand{\vd}{{\mathbf{d}}}

\newcommand{\vf}{{\mathbf{f}}}
\newcommand{\vg}{{\mathbf{g}}}

\newcommand{\vx}{{\mathbf{x}}}
\newcommand{\vy}{{\mathbf{y}}}
\newcommand{\vz}{{\mathbf{z}}}

\newcommand{\vA}{{\mathbf{A}}}
\newcommand{\vB}{{\mathbf{B}}}
\newcommand{\vC}{{\mathbf{C}}}

\newcommand{\vG}{{\mathbf{G}}}
\newcommand{\vH}{{\mathbf{H}}}
\newcommand{\vI}{{\mathbf{I}}}

\newcommand{\vM}{{\mathbf{M}}}

\newcommand{\vP}{{\mathbf{P}}}
\newcommand{\vQ}{{\mathbf{Q}}}

\newcommand{\vW}{{\mathbf{W}}}

\newcommand{\cE}{{\mathcal{E}}}

\newcommand{\cG}{{\mathcal{G}}}

\newcommand{\cV}{{\mathcal{V}}}

\newcommand{\RR}{\mathbb{R}}

\newcommand{\vzero}{\mathbf{0}}
\newcommand{\vone}{{\mathbf{1}}}

\newcommand{\st}{{\text{s.t.}}} 
\newcommand{\tr}{{\mathrm{tr}}} 

\newcommand{\Null}{\mathbf{Null}}
\newcommand{\Span}{\mathbf{span}}

\makeatletter
\let\@@span\span
\def\sp@n{\@@span\omit\advance\@multicnt\m@ne}
\makeatother

\DeclareMathOperator*{\Min}{minimize}

\DeclarePairedDelimiter{\dotp}{\langle}{\rangle}

\newcommand{\bc}{\begin{center}}
\newcommand{\ec}{\end{center}}

\newcommand{\bdm}{\begin{displaymath}}
\newcommand{\edm}{\end{displaymath}}

\newcommand{\beq}{\begin{equation}}
\newcommand{\eeq}{\end{equation}}

\newcommand{\bfl}{\begin{flushleft}}
\newcommand{\efl}{\end{flushleft}}

\newcommand{\bt}{\begin{tabbing}}
\newcommand{\et}{\end{tabbing}}

\newcommand{\beqn}{\begin{align}}
\newcommand{\eeqn}{\end{align}}

\newcommand{\beqs}{\begin{align*}} 
\newcommand{\eeqs}{\end{align*}}  

\newtheorem{assumption}{Assumption}

\newcommand{\Ker}{\mathbf{Ker}}
\newcommand{\Range}{\mathbf{Range}} 
\graphicspath{{figs/}}
\begin{document}

\title{On the linear convergence of two decentralized algorithms}


\author{Yao Li   \and  Ming Yan }

\institute{Yao Li \at
             Department of Mathematics, Michigan State University, 
              East Lansing, United States\\
              liyao6@msu.edu  
           \and
           Ming Yan  \at
              Department of Computational Mathematics, Science and Engineering and Department of Mathematics,              Michigan State University,               East Lansing, United States\\
              myan@msu.edu
}

\date{Received: date / Accepted: date}

\maketitle
\centerline{(Communicated by Zenon Mr\'oz)}

\begin{abstract}
Decentralized algorithms solve multi-agent problems over a connected network, where the information can only be exchanged with the accessible neighbors. Though there exist several decentralized optimization algorithms, there are still gaps in convergence conditions and rates between decentralized and centralized algorithms. In this paper, we fill some gaps by considering two decentralized algorithms: EXTRA and NIDS. They both converge linearly with strongly convex objective functions. We will answer two questions regarding them. What are the optimal upper bounds for their stepsizes? Do decentralized algorithms require more properties on the functions for linear convergence than centralized ones? More specifically, we relax the required conditions for linear convergence of both algorithms. For EXTRA, we show that the stepsize is comparable to that of centralized algorithms. For NIDS, the upper bound of the stepsize is shown to be exactly the same as the centralized ones. In addition, we relax the requirement for the objective functions and the mixing matrices. We provide the linear convergence results for both algorithms under the weakest conditions.
\end{abstract}


\section{Introduction}\label{sec: intro} 

We consider a network of $n$ agents minimizing $\bar{f}(x)\coloneqq\frac{1}{n}\sum_{i=1}^{n}f_{i}(x)$ collaboratively. 
Each convex and differentiable function $f_{i}\colon\RR^{p}\rightarrow\RR$ is known only by the corresponding agent $i$. The whole system is decentralized in the sense that each agent has an estimation of the global variable and can only exchange the estimation with their accessible neighbors during each iteration. A symmetric mixing matrix $\vW\in\RR^{n\times n}$ is used to encode the communication weights between the agents and enforce the consensus. The minimum condition for $\vW$ is having one eigenvalue $1$ with the all-one vector $\vone$ being a corresponding eigenvector. All other eigenvalues of $\vW$ are less than 1.

Early decentralized methods based on decentralized gradient descent~\cite{Nedic2009,ram2010distributed,nedic2010asynchronous,Kun2016,Jakovetic2014} have sublinear convergence for strongly convex objective functions, because of the diminishing stepsize that is needed to obtain a consensual and optimal solution. This sublinear convergence rate is much slower than that for centralized ones. The first decentralized algorithm with linear convergence~\cite{Shi2014} is based on Alternate Direction Multiplier Method (ADMM)~\cite{glowinski1975approximation,boyd2011distributed}. Note that this type of algorithms has $O(1/k)$ rate for general convex functions~\cite{Wei2013,Chang2015,Hong2015}. After that, many linearly convergent algorithms are proposed. Some examples are EXTRA~\cite{Shi2015}, NIDS~\cite{li2019decentralized}, DIGing~\cite{nedich2016achieving,Qu2016}, ESOM~\cite{mokhtari2016decentralized}, gradient tracking methods~\cite{zhu2010discrete,Xu2015,Lorenzo2016,Qu2016,nedich2016achieving,nedic2016geometrically,pu2018push}, exact diffusion~\cite{yuan2018exact1,yuan2018exact2}, and dual optimal~\cite{seaman2017optimal,uribe2018dual}. There are also works on composite functions, where each private function is the sum of a smooth and a nonsmooth functions~\cite{shi2015proximal,li2019decentralized,alghunaim2019linearly,Chen2012_2}. Another topic of interest is decentralized optimization over directed and dynamic graphs~\cite{Nedic2013,Xi2015,Zeng2015,sun2016distributed,nedich2016achieving,ling2013decentralized,Nedic2016stochastic}. Interested readers are referred to~\cite{Nedic2018} and the references therein.

This paper focuses on two linear convergent algorithms: EXTRA and NIDS, and provides better theoretical convergence results. EXact firsT-ordeR Algorithm (EXTRA) was proposed in~\cite{Shi2015}, and its iteration is described in~\eqref{extra}. For the general convex case, where each $f_i$ is convex and $L$-smooth (i.e., has a $L$-Lipschitz continuous gradient), the convergence condition in~\cite{Shi2015} is $\alpha \in \left(0,{1+\lambda_{\min}(\vW)\over L}\right)$. Therefore, there is an implicit condition for $\vW$ that the smallest eigenvalue of $\vW$ is larger than $-1$. Later the condition is relaxed to $\alpha \in \left(0,{5+3\lambda_{\min}(\vW)\over 4L}\right)$ in~\cite{li2017aprimal}, and the corresponding requirement for $\vW$ is that the smallest eigenvalue of $\vW$ is larger than $-5/3$. In addition, this condition for the stepsize is shown to be optimal, i.e., EXTRA may diverge if the condition is not satisfied. Though we can always manipulate $\vW$ to change the smallest eigenvalue, the convergence speed of EXTRA depends on the matrix $\vW$. In the numerical experiment, we will see that it is beneficial to choose small eigenvalues for EXTRA in certain scenarios.

The linear convergence of EXTRA requires additional conditions on the functions. There are mainly three types of conditions used in the literature: the strong convexity of $\bar f$ (and some weaker variants)~\cite{Shi2015}, the strong convexity of each $f_i$ (and some weaker variants)~\cite{li2017aprimal}, and the strong convexity of one function $f_i$~\cite{yuan2018exact2}. Note that the condition on $\bar f$ is much weaker than the other two; there are cases where $\bar f$ is strongly convex but none of $f_i$'s is. E.g., $f_{i}=\|e_{i}^{T}x\|^{2}_{2}$ for $p=n>1$, where $e_{i}$ is the vector whose $i$th component is $1$ and all other components are $0$. If $\bar f$ is (restricted) strongly convex with parameter $\mu_{\bar f}$, the linear convergence of EXTRA is shown when $\alpha \in (0,{\mu_{\bar f}(1+\lambda_{\min}(\vW))\over L^2})$ in~\cite{Shi2015}. The upper bound for the stepsize is very conservative, and the better performance with a larger stepsize was shown numerically in~\cite{Shi2015} without proof. If each $f_i$ is strongly convex with parameter $\mu$, the linear convergence is shown when $\alpha \in \left(0,{1+\lambda_{\min}(\vW)\over L+\mu}\right)$ and $\alpha \in \left(0,{5+3\lambda_{\min}(\vW)\over 4L}\right)$ in~\cite{alghunaim2019linearly} and~\cite{li2017aprimal}, respectively. One contribution of this paper is to show the linear convergence of EXTRA under the (restricted) strong convexity of $\bar f$ and $\alpha \in \left(0,{5+3\lambda_{\min}(\vW)\over 4L}\right)$.

The algorithm NIDS (Network InDependent Stepsize) was proposed in~\cite{li2019decentralized}. Though there is a small difference from EXTRA, NIDS can choose a stepsize that does not depend on the mixing matrices. The linear convergence of NIDS in~\cite{li2019decentralized} requires $\vI\succcurlyeq \vW \succ -\vI$ and strong convexity of $\vf(\vx)$. In this paper, we relax this condition for linear convergence to (restricted) strong convexity of $\bar f(x)$ and the relaxed mixing matrices with $\vI\succcurlyeq \vW \succ -(5/3)\vI$.

In sum, we provide new and stronger linear convergence results for two state-of-the-art algorithms: EXTRA and NIDS. More specifically,
\begin{itemize}
\item We show the linear convergence of EXTRA under the strong convexity of $\bar f$ and the relaxed condition $\lambda_{\min}(\vW)>-5/3$.
The upper bound of the stepsize can be as large as ${5+3\lambda_{\min}(\vW)\over 4L}$, which is shown to be optimal in~\cite{li2017aprimal} for general convex problems;
\item We show the linear convergence of NIDS under the same condition on $\bar f$ and $\vW$ as EXTRA with any network-independent stepsize $\alpha\in (0,2/L)$.
\end{itemize}


\section{Notation}\label{sec:notation}
We let
\beq
\vf(\vx)\coloneqq\sum_{i=1}^{n}f_{i}(x_{i}), \label{pbm}
\eeq
where each $x_{i}\in\RR^{p}$ is the local copy of the global variable $x$ and the $k$th iterated point is $x_{i}^{k}$. Since agent $i$ has its own estimate $x_i$ of the global variable $x$, we put them together and define
$$\vx = [x_1,x_2,\cdots,x_n]^\top\in\RR^{n\times p}.$$
The gradient of $\vf$ is defined as
$$\nabla \vf(\vx) = [\nabla f_1(x_1),\nabla f_2(x_2),\cdots, \nabla f_n(x_n)]^\top\in \RR^{n\times p}.$$
We say that $\vx$ is consensual if $x_{1}=x_{2}=\cdots =x_{n},$ i.e., $\vx=\vone x^{\top}$, where $x\in\RR^{p\times 1}$ and $\vone=[1,1,\cdots,1]^{\top}\in\RR^{n\times 1}.$

In this paper, we use $\|\cdot\|$ and $\dotp{\cdot,\cdot}$ to denote the Frobenious norm and the corresponding inner product, respectively.
For a given matrix $\vM\in \RR^{n\times p}$ and any positive (semi)definite matrix $\vH$, which is denoted as $\vH\succ 0$ ($\vH\succcurlyeq 0$ for positive semidefinite), we define $\|\vM\|_{\vH}\coloneqq\sqrt{\tr(\vM^{\top}\vH\vM)}.$
The largest and the smallest eigenvalues of a matrix $\vA$ are defined as $\lambda_{\max}(\vA)$ and $\lambda_{\min}(\vA)$.
For a symmetric positive semidefinite matrix $\vA$, we let $\lambda_{\min}^+(\vA)$ be the smallest nonzero eigenvalue. $\vA^{\dagger}$ is the pseudo inverse of $\vA$. For a matrix $\vA\in\RR^{n\times n},$ we say a matrix $\vB\in\RR^{n\times p}$ is in $\Ker\{\vA\}$ if $\vA\vB=\vzero_{n\times p},$ and $\vB$ is in $\Range\{\vA\}$ if there exists $\vC\in\RR^{n\times p}$ such that $\vB=\vA\vC.$
For simplicity, we may use $\vx^+$ and $\vx$ to replace $\vx^{k+1}$ and $\vx^k$, respectively, in the proofs.

\section{Algorithms and Prerequisites}

One iteration of EXTRA can be expressed as
\beq\label{extra}
\begin{aligned}
\vx^{k+2}= (\vI+\vW)\vx^{k+1}-\widetilde{\vW}\vx^{k} -\alpha[\nabla\vf(\vx^{k+1})-\nabla\vf(\vx^{k})].
\end{aligned}
\eeq
The stepsize $\alpha>0$, and the symmetric matrices $\vW$ and $\widetilde\vW$ satisfy $\vI+\vW\succcurlyeq 2\widetilde\vW\succcurlyeq 2\vW$.
The initial value $\vx^{0}$ is chosen arbitrarily, and $\vx^{1}=\vW\vx^{0}-\alpha\nabla\vf(\vx^{0})$.
In practice, we usually let $\widetilde{\vW}={\vI+\vW\over 2}.$

One iteration of NIDS is
\beq\label{nids}
\begin{aligned}
\textstyle\vx^{k+2}= \frac{\vI+\vW}{2}\left[2\vx^{k+1}-\vx^{k}\right. \left.-\alpha(\nabla\vf(\vx^{k+1})-\nabla\vf(\vx^{k}))\right],
\end{aligned}
\eeq
where $\alpha>0$ is the stepsize.
The initial value $\vx^{0}$ is chosen arbitrarily, and $\vx^{1}={\vI+\vW\over 2}[\vx^{0}-\alpha\nabla\vf(\vx^{0})]$.

If we choose $\widetilde{\vW}=\frac{\vI+\vW}{2}$ in~\eqref{extra}, the difference between EXTRA and NIDS in the above mathematical forms happens only in the communicated data, i.e.,whether we exchange the gradient information or not at each step. In practice, EXTRA can gain the advantage of time overlap by parallelizing communication and gradient evaluation, while NIDS evaluates the gradient and then communicates after the gradient is added.
However, this small difference brings big changes in the convergence~\cite{li2019decentralized}.
In order for both algorithms to converge, we have the following assumptions on $\vW$ and $\widetilde\vW$.

\begin{assumption}[Mixing matrix]\label{asm1}
The connected network $\cG=\{\cV,\cE\}$ consists of a set of nodes $\cV=\{1,2,\cdots,n\}$ and a set of undirected edges $\cE$.
An undirected edge $(i,j)\in\cE$ means that there is a connection between agents $i$ and $j$ and both agents can exchange data.
The mixing matrices $\vW=[w_{ij}]\in \RR^{n\times n}$ and $\widetilde{\vW}=[\widetilde{w}_{ij}]\in \RR^{n\times n}$ satisfy:
\begin{enumerate}
\item  \text{(Decentralized property):}\ If $i\neq j$ and $(i,j)\notin \cE$, then $w_{ij}=\widetilde{w}_{ij}=0$\label{part1}.
\item  \text{(Symmetry):}\ $\vW = \vW^\top$, $\widetilde{\vW}=\widetilde{\vW}^{\top}$\label{part2}.
\item  \text{(Null space property):} $\Null\{\vW-\widetilde{\vW}\}=\Span\{\vone\} \subseteq \Null\{\vI-\widetilde{\vW}\}.$
    \label{part3}
\item  \text{(Spectral property):} $\frac{\vI+\vW}{2}\succcurlyeq\widetilde{\vW}\succ-\frac{1}{3}\vI,\quad \widetilde{\vW}\succcurlyeq\vW.$
\label{part4}
\end{enumerate}
\end{assumption}
\begin{remark}\label{rm1}
Parts~\ref{part2}-\ref{part4} imply that the spectrum of $\vW$ is enlarged to $(-\frac{5}{3},1]$, while the original assumption is $(-1,1]$ for doubly stochastic matrices.
Therefore, in our assumption, $\frac{\vI+\vW}{2}$ does not have to be positive definite.
This assumption for $\vW$ is strictly weaker than those in~\cite{Shi2015} and~\cite{li2019decentralized}.
\end{remark}
\begin{remark}\label{rmW}
   From~\cite[Proposition 2.2]{Shi2015}, $\Null\{\vI-\vW\}=\Span\{\vone\}.$ It is a critical result for both algorithms.
\end{remark}

Before showing their theoretical results, we reformulate both algorithms.

\noindent{\bf Reformulation of EXTRA:}
We reformulate EXTRA by introducing a variable $\vy\in\RR^{n\times p}$ as
\begin{subequations}\label{reformulate}
\begin{align}
&\vx^{k+1}=\widetilde{\vW}\vx^{k}+\vy^{k}-\alpha\nabla\vf(\vx^{k})\label{reformulate:a},\\
&\vy^{k+1}=\vy^{k}-(\widetilde{\vW}-\vW)\vx^{k+1}\label{reformulate:b},
\end{align}
\end{subequations}
with $\vy^{0}=-(\widetilde{\vW}-\vW)\vx^{0}$. Then~\eqref{reformulate} is equivalent to EXTRA~\eqref{extra}.

\begin{proposition}\label{pp1}
Let the $\vx$-sequence generated by~\eqref{reformulate} with $\vy^0=-(\widetilde{\vW}-\vW)\vx^{0}$ be $\{\vx^{k}\}_{k=1}^{\infty}$, then it's identical to the sequence generated by EXTRA~\eqref{extra} with the same initial point $\vx^{0}$. 
\end{proposition}
\begin{proof}
From~\eqref{reformulate:a}, we have
\begin{align*}
\vx^1= &\widetilde{\vW}\vx^{0}+\vy^{0}-\alpha\nabla\vf(\vx^{0})= \widetilde{\vW}\vx^{0}-(\widetilde{\vW}-\vW)\vx^{0}-\alpha\nabla\vf(\vx^{0}) \\
= &\vW\vx^{0}-\alpha\nabla\vf(\vx^{0}).
\end{align*}
For $k\geq 0$, we have
\begin{equation*}
\begin{aligned}
  \vx^{k+2}= &\widetilde{\vW}\vx^{k+1}+\vy^{k+1}-\alpha\nabla\vf(\vx^{k+1})           = \vW\vx^{k+1}+\vy^{k}-\alpha\nabla\vf(\vx^{k+1})\\
          = &(\vI+\vW)\vx^{k+1}-\widetilde{\vW}\vx^{k} -\alpha[\vf(\vx^{k+1})-\vf(\vx^{k})],
  \end{aligned}
  \end{equation*}
	where the second and last equalities are from~\eqref{reformulate:b} and~\eqref{reformulate:a}, respectively.  \qed
\end{proof}

\begin{remark}\label{rm2}
  By~\eqref{reformulate:b} and the assumption of $\vy^{0}$, $\vy^{k}\in \Range\{\widetilde{\vW}-\vW\}$ for all $k$. Also, $\vx^{k+1}=(\widetilde{\vW}-\vW)^{\dagger}(\vy^{k}-\vy^{k+1})+\vz^{k+1}$ for some $\vz^{k+1}\in\Ker\{\widetilde{\vW}-\vW\}$.
\end{remark}

\noindent{\bf Reformulation of NIDS:}
We adopt the following reformulation from~\cite{li2019decentralized}:
\begin{subequations}\label{reform_nids}
\begin{align}
&\textstyle \vd^{k+1}=\vd^{k}+\frac{\vI-\vW}{2\alpha}[\vx^{k}-\alpha\nabla\vf(\vx^{k})-\alpha\vd^{k}]\label{reform_nids:a},\\
&\vx^{k+1}=\vx^{k}-\alpha\nabla\vf(\vx^{k})-\alpha\vd^{k+1}\label{reform_nids:b},
\end{align}
\end{subequations}
with $\vd^0=\vzero.$
The equivalence is shown in~\cite{li2019decentralized}.

To establish the linear convergence of EXTRA and NIDS, we need the following two assumptions.
\begin{assumption}[Uniqueness]\label{asm3}
There is a unique minimizer $x^*$ for $\bar f(x)$.
\end{assumption}

\begin{assumption}[$L$-smoothness and restricted strong convexity]\label{LipResSC}
Each function $f_{i}$ is a proper, closed and convex function with a Lipschitz continuous gradient:
\beq\label{Lip1}
  \|\nabla f_{i}(x)-\nabla f_{i}(\widetilde{x})\|\leq L\|x-\widetilde{x}\|,~\forall x,~\widetilde{x}\in\RR^{p},
\eeq
where $L>0$ is the Lipschitz constant.
Furthermore, $\bar{f}(x)$ is (restricted) strongly convex with respect to $x^{*}$:
\beq\label{Stc1}
  \dotp{x-x^{*},\nabla \bar{f}(x)-\nabla \bar{f}(x^{*})}\geq\mu_{\bar{f}}\|x-x^{*}\|^{2},~\forall x\in\RR^{p}.
\eeq
\end{assumption}

From~\cite[Theorem 2.1.5]{nesterov2013introductory}, the inequality~\eqref{Lip1} is equivalent to, for any $\vx,~\tilde{\vx}\in\RR^{n\times p}$,
\begin{equation}\label{Lip3}
  \dotp{\vx-\widetilde{\vx},\nabla \vf(\vx)-\nabla \vf(\widetilde{\vx})}\geq L^{-1}\|\nabla\vf(\vx)-\nabla\vf(\widetilde{\vx})\|^{2}.
\end{equation}

\begin{proposition}[{\cite[Appendix A]{Shi2015}}]\label{prop2}
The following two statements are equivalent:
\begin{enumerate}
  \item $\bar{f}(x)$ is (restricted) strongly convex with respect to $x^{*};$
  \item For any $\eta>0$, $\vg(\vx)\coloneqq\vf(\vx)+\frac{\eta}{2}\|\vx\|^{2}_{\vI-\vW}$ is $\mu_{\vg}$ (restricted) strongly convex with respect to $\vx^{*}=\vone(x^{*})^{\top}$. Specially, we can characterize
	$$ \mu_{\vg}=\min\left\{\frac{\mu_{\bar{f}}}{2}, \frac{\mu_{\bar{f}}^{2}\lambda^+_{\min}(\vI-\vW)}{\mu_{\bar{f}}^{2}+16L^{2}}\eta\right\}.$$
\end{enumerate}
\end{proposition}

This proposition shows
\begin{align}
   \dotp{\vx-\vx^{*},\nabla \vf(\vx)-\nabla \vf(\vx^{*})}+\eta\|\vx-\vx^{*}\|^{2}_{\vI-\vW}\geq  \mu_{\vg}\|\vx-\vx^{*}\|^{2}\label{Stc2}
\end{align}
for any $\vx\in\RR^{n\times p}.$

\section{New Linear Convergence Results for EXTRA and NIDS}
Throughout this section, we assume that Assumptions~\ref{asm1}-\ref{LipResSC} hold. Two techniques are used to show the linear convergence: a) Proposition~\ref{prop2} serves as a bridge to connect $\bar{f}(x)$ and $\vf(\vx)$. It is the key to the weaker assumption on objective functions. b) Both algorithms are equivalent to the extended Proximal Alternating Predictor-Corrector (PAPC) in~\cite{li2017aprimal}, and this equivalence is the key to relaxing the conditions on the mixing matrices $\vW$ and $\widetilde\vW$.

\subsection{Linear Convergence of EXTRA}
When $\widetilde{\vW}=\frac{\vI+\vW}{2}$, EXTRA is recovered by applying the extended PAPC in~\cite{li2017aprimal} to the following dual form of the decentralized consensus problem $$\Min_\vy~ \vf^*(\sqrt{\vI-\vW}\vy),$$ 
where $\vf^*$ is the conjugate function of $\vf$ and $\vy$ is the dual variable. In this case, EXTRA has the optimal bound of the stepsize over the relaxed mixing matrix $\vW\succ-(5/3)\vI$. This fact enlightens us on the critical Lemma~\ref{keyineqlem}.

For simplicity, we introduce some notations.
Because of part~\ref{part4} of Assumption~\ref{asm1}, given the mixing matrices $\vW$ and $\widetilde{\vW},$ there is a constant
$$\theta\in\Big(\frac{3}{4},\min\big\{\frac{1}{1-\lambda_{\min}(\widetilde\vW)},1\big\}\Big]$$ such that
\begin{align}
\overbar\vW  \coloneqq & \theta\widetilde{\vW}+(1-\theta)\vI\succ\vzero,\label{barW}\\
\vH\coloneqq & \textstyle \overbar{\vW}+(\theta-\frac{1}{2})(\vI-\widetilde{\vW})=\frac{\vI+\widetilde{\vW}}{2}\succ\vzero,\\
\vM\coloneqq & (\widetilde{\vW}-\vW)^{\dagger}\succcurlyeq\vzero,\label{eq:defM}\\
\vG\coloneqq & \vW+\vI-2\widetilde{\vW}\succcurlyeq\vzero.
\end{align}
Based on~\eqref{barW}, we have
\beq\label{barWeq}
\widetilde{\vW}=\overbar{\vW}-(1-\theta)(\vI-\widetilde{\vW}).
\eeq
Let $(\vx^{*},\vy^{*})$ be a fixed point of~\eqref{reformulate}, it is straightforward to show that 
  \begin{align}
  (\widetilde{\vW}-\vW)\vx^{*}= & \vzero\label{soln1:a}. 
  \end{align}
Part~\ref{part3} of Assumption~\ref{asm1} shows that $\vx^*$ is consensual, i.e., $\vx^{*}=\vone(x^{*})^{\top}$ for certain $x^{*}\in\mathbb{R}^{p}$.
The $\vy$-iteration in~\eqref{reformulate:b} and the initialization of $\vy^0$ show $\vy^{k}\in \Range\{\widetilde{\vW}-\vW\}=\Ker\{\vone^{\top}\}$.
Then we have $\vone^{\top}\vy^{*}=\alpha\vone^{\top}\nabla\vf(\vx^{*})=0$.
Thus, $x^{*}$ is the unique minimizer of $\bar f(x)$.

\begin{lemma}[Norm over range space {\cite[Lemma 3]{li2019decentralized}}]\label{lem:norm over range space}
For any symmetric positive (semi)definite matrix $\vA\in \RR^{n\times n}$ with rank $r$ ($r\leq n$),
let $\lambda_1\geq \lambda_2\geq \dots\geq \lambda_r>0$ be its $r$ eigenvalues.
Then $\Range\{\vA\}$ is a $rp$-dimensional subspace in $\RR^{n\times p}$ and has a norm defined by $\|\vx\|^2_{\vA^{\dagger}}:=\langle \vx, \vA^{\dagger} \vx\rangle$, where $\vA^{\dagger}$ is the pseudo inverse of $\vA$.
In addition, $\lambda_1^{-1}\|\vx\|^2\leq \|\vx\|^2_{\vA^{\dagger}}\leq \lambda_r^{-1}\|\vx\|^2$ for all $\vx\in\Range\{\vA\}$.
\end{lemma}

For simplicity, we let $\vx^+$ and $\vx$ stand for $\vx^{k+1}$ and $\vx^{k}$, respectively, in the proofs.
The same simplification applies to $\vy^k$.
\begin{lemma}[Norm equality]\label{lem:norm equality}
Let $\{(\vx^{k},\vy^{k})\}_{k=1}^{\infty}$ be the sequence generated by~\eqref{reformulate}, then it satisfies
\beq\label{ne}
\|\vx^{k+1}-\vx^{*}\|^{2}_{\widetilde{\vW}-\vW}=\|\vy^{k}-\vy^{k+1}\|^{2}_{\vM}.
\eeq
\end{lemma}
\begin{proof}
From Remark~\ref{rm2}, we have
\begin{equation}
\vx^{+}=\vM(\vy-\vy^{+})+\vz^{+} \label{eq:lemma2eq1}
\end{equation}
for $\vz^+\in\Ker\{\widetilde{\vW}-\vW\}.$
This equality and~\eqref{soln1:a} give
\begin{equation*}
\begin{aligned}
\|\vx^{+}-\vx^{*}\|^{2}_{\widetilde{\vW}-\vW}
 =& \dotp{\vx^{+}-\vx^{*},(\widetilde{\vW}-\vW)(\vx^{+}-\vx^{*})}= \dotp{\vx^{+},(\widetilde{\vW}-\vW)\vx^{+}} \\
 = & \dotp{\vM(\vy-\vy^{+}),\vy-\vy^{+}}= \|\vy-\vy^{+}\|_{\vM}^{2},
\end{aligned}
\end{equation*}
where the third equality holds because of~\eqref{eq:defM},~\eqref{eq:lemma2eq1}, and $\vy-\vy^{+}\in \Range(\widetilde\vW-\vW)$.
\qed
\end{proof}




\begin{lemma}[A key inequality for EXTRA]\label{keyineqlem}
Let $\{(\vx^{k},\vy^{k})\}_{k=1}^{\infty}$ be the sequence generated by~\eqref{reformulate}, then we have
  \begin{align}\label{keyineqfma}
     & \|\vx^{k+1}-\vx^{*}\|^{2}_{\vH}+\|\vy^{k+1}-\vy^{*}\|_{\vM}^{2}\nonumber \\
 \leq& \|\vx^{k}-\vx^{*}\|^{2}_{\vH}+\|\vy^{k}-\vy^{*}\|^{2}_{\vM} -\|\vx^{k}-\vx^{k+1}\|^{2}_{(\theta-\frac{3}{4})(\vI-\widetilde{\vW})} -\|\vx^{k+1}-\vx^{*}\|^{2}_{\vG}\nonumber\\
     & -\|\vx^{k}-\vx^{k+1}\|^{2}_{\overbar{\vW}}-2\alpha\dotp{\vx^{k+1}-\vx^{*},\nabla\vf(\vx^{k})-\nabla\vf(\vx^{*})}.
  \end{align}
\end{lemma}
\begin{proof} The iteration~\eqref{reformulate} and equation~\eqref{barWeq} show
\begin{align}
   &2\alpha\dotp{\vx^{+}-\vx^{*},\nabla\vf(\vx)-\nabla\vf(\vx^{*})}\nonumber\\
  =& 2\dotp{\vx^{+}-\vx^{*}, \widetilde{\vW}(\vx-\vx^{+})+\widetilde{\vW}(\vx^{+}-\vx^{*}) -(\vx^{+}-\vx^{*})+(\vy-\vy^{*})}\nonumber\\
  =& 2\dotp{\vx^{+}-\vx^{*}, \widetilde{\vW}(\vx-\vx^{+}) +(\widetilde{\vW}-\vI)(\vx^{+}-\vx^{*}) \nonumber\\
   & +(\widetilde{\vW}-\vW)(\vx^{+}-\vx^{*})+\vy^{+}-\vy+\vy-\vy^{*}}\nonumber\\
  =& 2\dotp{\vx^{+}-\vx^{*},\widetilde{\vW}(\vx-\vx^{+})} +2\dotp{\vx^{+}-\vx^{*},\vy^{+}-\vy^{*}} -2\|\vx^{+}-\vx^{*}\|^{2}_{\vG}\nonumber\\
  =& 2\dotp{\vx^{+}-\vx^{*},\overbar{\vW}(\vx-\vx^{+})}-2\dotp{\vx^{+}-\vx^{*},(1-\theta)(\vI-\widetilde{\vW})(\vx-\vx^{+})}\nonumber\\
   & +2\dotp{\vx^{+}-\vx^{*},\vy^{+}-\vy^{*}}-2\|\vx^{+}-\vx^{*}\|^{2}_{\vG}, \label{eq:Lemma3eq1}
\end{align}
where the first equality comes from~\eqref{reformulate:a}, the second one follows~\eqref{reformulate:b}, and the last one is from~\eqref{barWeq}. From Remark~\ref{rm2}, $\vx^{+}-\vx^{*}=\vM(\vy-\vy^{+})+\vz^{+}-\vx^{*}$ for some $\vz^{+}\in\Ker\{\widetilde{\vW}-\vW\}$.
Thus $\dotp{\vz^{+}-\vx^{*},\vy^{+}-\vy^{*}}=0,$ and the equality~\eqref{eq:Lemma3eq1} can be rewritten as
  \begin{equation*}
    \begin{aligned}
      & 2\alpha\dotp{\vx^{+}-\vx^{*},\nabla\vf(\vx)-\nabla\vf(\vx^{*})}\\
     =& 2\dotp{\vx^{+}-\vx^{*},\overbar{\vW}(\vx-\vx^{+})}  -2\dotp{\vx^{+}-\vx^{*},(1-\theta)(\vI-\widetilde{\vW})(\vx-\vx^{+})}\\
      & +2\dotp{\vM(\vy-\vy^{+}),\vy^{+}-\vy^{*}}  -2\|\vx^{+}-\vx^{*}\|^{2}_{\vG}.
    \end{aligned}
  \end{equation*}

  Using the basic equality
 $
  2\dotp{a-b,b-c}
  =\|a-c\|^{2}-\|a-b\|^{2}-\|b-c\|^{2} $
  and Lemma~\ref{lem:norm equality}, we have
    \begin{align}
     & \|\vx^{+}-\vx^{*}\|^{2}_{\overbar{\vW}}-\|\vx^{+}-\vx^{*}\|^{2}_{(1-\theta)(\vI-\widetilde{\vW})} +\|\vy^{+}-\vy^{*}\|^{2}_{\vM}\nonumber\\
    =& \|\vx-\vx^{*}\|^{2}_{\overbar{\vW}}-\|\vx-\vx^{*}\|^{2}_{(1-\theta)(\vI-\widetilde{\vW})}+\|\vy-\vy^{*}\|^{2}_{\vM}\nonumber\\
    &-\|\vx-\vx^{+}\|^{2}_{\overbar{\vW}} +\|\vx-\vx^{+}\|^{2}_{(1-\theta)(\vI-\widetilde{\vW})}-\|\vx^{+}-\vx^{*}\|^{2}_{\widetilde{\vW}-\vW}\nonumber\\
    &-2\|\vx^{+}-\vx^{*}\|^{2}_{\vG}-2\alpha\dotp{\vx^{+}-\vx^{*},\nabla\vf(\vx)-\nabla\vf(\vx^{*})}.\label{eq1}
    \end{align}
  Note that the following inequality holds,
  \begin{equation*}
    \begin{aligned}
 \textstyle \frac{1}{2}\|\vx^{+}-\vx^{*}\|^{2}_{\widetilde{\vW}-\vW}
    \leq \|\vx^{+}-\vx^{*}\|^{2}_{\widetilde{\vW}-\vW} 	+\frac{1}{2}\|\vx-\vx^{*}\|^{2}_{\widetilde{\vW}-\vW}
         -\frac{1}{4}\|\vx-\vx^{+}\|^{2}_{\widetilde{\vW}-\vW}.
    \end{aligned}
  \end{equation*}
Adding it onto both sides of~\eqref{eq1}, we have
    \begin{align}
     & \textstyle\|\vx^{+}-\vx^{*}\|^{2}_{\vH}-\frac{1}{2}\|\vx^{+}-\vx^{*}\|^{2}_{\vG} +\|\vy^{+}-\vy^{*}\|^{2}_{\vM}\nonumber\\
    \leq& \textstyle\|\vx-\vx^{*}\|^{2}_{\vH}-\frac{1}{2}\|\vx-\vx^{*}\|^{2}_{\vG}+\|\vy-\vy^{*}\|^{2}_{\vM}\nonumber\\
     & -\|\vx-\vx^{+}\|^{2}_{\overbar{\vW}}-\|\vx-\vx^{+}\|^{2}_{(\theta-\frac{3}{4})(\vI-\widetilde{\vW})} +\frac{1}{4}\|\vx-\vx^{+}\|^{2}_{\vG}\nonumber\\
     &-2\|\vx^{+}-\vx^{*}\|^{2}_{\vG} -2\alpha\dotp{\vx^{+}-\vx^{*},\nabla\vf(\vx)-\nabla\vf(\vx^{*})}.\label{eq2}
    \end{align}
Apply the inequality
  $\textstyle\frac{1}{4}\|\vx-\vx^{+}\|^{2}_{\vG}\leq\frac{1}{2}\|\vx-\vx^{*}\|^{2}_{\vG}+\frac{1}{2}\|\vx^{+}-\vx^{*}\|^{2}_{\vG},$
  then the key inequality~\eqref{keyineqfma} is obtained. \qed
\end{proof}

In the following theorem, we assume $\vG\neq\vzero$ (i.e., $\widetilde{\vW}\not=({\vI+\vW})/{2}$).
It is easy to amend the proof to show the result for this special case.

  \begin{theorem}[Q-linear convergence of EXTRA]
    Under Assumptions~\ref{asm1}-\ref{LipResSC}, we define
		\begin{align}
      r_{1}= & \textstyle \frac{4\theta-3}{4(1-\theta)^{2}\lambda_{\max}(\overbar{\vW}^{-1}(\vI-\widetilde{\vW}))}>0, \label{r1}\\
      r_{2}= & \textstyle \frac{1}{2\lambda_{\max}(\vG\overbar{\vW}^{-1})}>0, \label{r2} \\
      r_{3}= & \textstyle \frac{r_1r_2}{r_{1}+r_{2}+r_{1}r_{2}}\in(0,1),   \label{r3}
    \end{align}
		and choose two small parameters $\xi$ and $\eta$ such that
    \begin{align}
		\xi\in  & \textstyle \Big(0,\min\Big\{\frac{r_{3}}{4\lambda_{\max}(\overbar{\vW}\vM)},1\Big\}\Big), \label{eps}\\
    \eta\in & \textstyle \Big(0,{\lambda_{\min}(\overbar\vW)\xi\over 4\alpha\lambda_{\min}(\overbar\vW)-2\alpha^2L}\Big).\label{eta}
    \end{align}
			In addition, we define
		\begin{align*}
		\vP\coloneqq & \textstyle \vH+\frac{\xi}{2}(\vI-\vW)\succ\vzero,\\
		\vQ\coloneqq & \vM+(r_{3}-2\xi\lambda_{\max}(\overbar{\vW}\vM))\overbar{\vW}^{-1}\succ\vzero.
		\end{align*}
Then for any stepsize $\alpha\in(0,\frac{2\lambda_{\min}(\overbar{\vW})}{L}),$ we have
    \begin{equation}\label{linineq}
    \begin{aligned}
       \|\vx^{k+1}-\vx^{*}\|_{\vP}^{2}+\|\vy^{k+1}-\vy^{*}\|^{2}_{\vQ}\leq \rho(\|\vx^{k}-\vx^{*}\|_{\vP}^{2}+\|\vy^{k}-\vy^{*}\|^{2}_{\vQ}),
    \end{aligned}
    \end{equation}
		where
    \begin{equation}\label{rho}
    \begin{aligned}
    \rho\coloneqq\max& \textstyle \Big\{\ 1-\Big(2\alpha-\frac{\alpha^{2}L}{\lambda_{\min}(\overbar{\vW})}\Big){\mu_{\vg}},  \Big(4\alpha-\frac{2\alpha^{2}L}{\lambda_{\min}(\overbar{\vW})}\Big)\frac{\eta}{\xi}, \\
     &\textstyle 1-\frac{r_{3}-4\xi\lambda_{\max}(\overbar{\vW}\vM)}{r_{3}+(1-2\xi)\lambda_{\max}(\overbar{\vW}\vM)}\Big\}<1.
    \end{aligned}
    \end{equation}
  \end{theorem}
\begin{proof}
From~\eqref{keyineqfma} in Lemma~\ref{keyineqlem}, we have
\begin{align}
		& \|\vx^{+}-\vx^{*}\|^{2}_{\vH}+\|\vy^{+}-\vy^{*}\|_{\vM}^{2}\nonumber\\
\leq& \|\vx-\vx^{*}\|^{2}_{\vH}+\|\vy-\vy^{*}\|^{2}_{\vM} -\|\vx-\vx^{+}\|^{2}_{(\theta-\frac{3}{4})(\vI-\widetilde{\vW})}-\|\vx^{+}-\vx^{*}\|^{2}_{\vG}\nonumber\\
		&-\|\vx-\vx^{+}\|^{2}_{\overbar{\vW}} -2\alpha\dotp{\vx^{+}-\vx^{*},\nabla\vf(\vx)-\nabla\vf(\vx^{*})}.\label{thmeq1}
\end{align}
Then we find an upper bound of $-\|\vx-\vx^{+}\|_{\overbar{\vW}}^{2}-2\alpha\langle\vx^{+}-\vx^{*}, \nabla\vf(\vx)-\nabla\vf(\vx^{*})\rangle$.
\begin{equation*}\label{thmeq2}
	\begin{aligned}
		& -\|\vx-\vx^{+}\|^{2}_{\overbar{\vW}} -2\alpha\langle\vx^+-\vx^{*}, \nabla\vf(\vx)-\nabla\vf(\vx^{*})\rangle\\
	 =& \alpha^{2}\|\nabla\vf(\vx)-\nabla\vf(\vx^{*})\|^{2}_{\overbar{\vW}^{-1}}-2\alpha\dotp{\vx-\vx^*,\nabla\vf(\vx)-\nabla\vf(\vx^{*})}\\
		& -\|\overbar{\vW}(\vx-\vx^+)-\alpha(\nabla\vf(\vx)-\nabla\vf(\vx^{*}))\|^{2}_{\overbar{\vW}^{-1}}\\
\leq& \textstyle -\big(2\alpha-\frac{\alpha^{2}L}{\lambda_{\min}(\overbar{\vW})}\big)\dotp{\vx-\vx^{*},\nabla\vf(\vx)-\nabla\vf(\vx^{*})}\\
		& -\|\overbar{\vW}(\vx-\vx^+)-\alpha(\nabla\vf(\vx)-\nabla\vf(\vx^{*}))\|^{2}_{\overbar{\vW}^{-1}},
	\end{aligned}
\end{equation*}
where, the inequality comes from~\eqref{Lip3}. 
Combining it with~\eqref{thmeq1}, we have
\begin{align}
	& \|\vx^{+}-\vx^{*}\|^{2}_{\vH}+\|\vy^{+}-\vy^{*}\|_{\vM}^{2} -\|\vx-\vx^{*}\|^{2}_{\vH}-\|\vy-\vy^{*}\|^{2}_{\vM}\nonumber\\
\leq& \textstyle- \big(2\alpha-\frac{\alpha^{2}L}{\lambda_{\min}(\overbar{\vW})}\big)\dotp{\vx-\vx^{*},\nabla\vf(\vx)-\nabla\vf(\vx^{*})} \nonumber\\
	& -\|\overbar{\vW}(\vx-\vx^+)-\alpha(\nabla\vf(\vx)-\nabla\vf(\vx^{*}))\|^{2}_{\overbar{\vW}^{-1}}\nonumber\\
	& -\|\vx-\vx^{+}\|^{2}_{(\theta-\frac{3}{4})(\vI-\widetilde{\vW})}-\|\vx^{+}-\vx^{*}\|^{2}_{\vG}.\label{thmeq3}
\end{align}
The inequality~\eqref{thmeq3} shows that $\{(\vx^{k},\vy^{k})\}_{k=1}^{\infty}$ is a Cauchy sequence converging to the fixed point $(\vx^{*},\vy^{*})$ of~\eqref{reformulate}.
From~\eqref{Stc2}, we can bound the first term on the right hand side of~\eqref{thmeq3} as
\begin{align}
    &\textstyle -\big(2\alpha-\frac{\alpha^{2}L}{\lambda_{\min}(\overbar{\vW})}\big)\dotp{\vx-\vx^{*},\nabla\vf(\vx)-\nabla\vf(\vx^{*})}\nonumber\\
\leq&\textstyle \big(2\alpha-\frac{\alpha^{2}L}{\lambda_{\min}(\overbar{\vW})}\big)\eta\|\vx-\vx^*\|_{\vI-\vW}^2-\big(2\alpha-\frac{\alpha^{2}L}{\lambda_{\min}(\overbar{\vW})}\big)\mu_g\|\vx-\vx^*\|^2.\label{eq:ub1}
\end{align}
Next, we bound the two terms involving successive iterated points, i.e.,  $-\|\overbar{\vW}(\vx-\vx^+)-\alpha(\nabla\vf(\vx)-\nabla\vf(\vx^{*}))\|^{2}_{\overbar{\vW}^{-1}}-\|\vx-\vx^{+}\|^{2}_{(\theta-\frac{3}{4})(\vI-\widetilde{\vW})}$.
Note that
\begin{align}
		& \overbar{\vW}(\vx-\vx^+)-\alpha(\nabla\vf(\vx)-\nabla\vf(\vx^{*}))\nonumber\\
	= & \vG(\vx^{+}-\vx^{*})-(\vy^{+}-\vy^{*}) +(1-\theta)(\vI-\widetilde{\vW})(\vx-\vx^+). \label{sucterm}
\end{align}
We use $T_{1},~T_{2}$, and $T_{3}$ to denote the three terms on the right hand side of~\eqref{sucterm}, respectively.
Using the definition of $r_{1}$ in~\eqref{r1}, we have
\begin{equation*}\label{thmeq4}
	\begin{aligned}
		& -\|T_{1}+T_{2}+T_{3}\|^{2}_{\overbar{\vW}^{-1}}-\|\vx-\vx^{+}\|^{2}_{(\theta-\frac{3}{4})(\vI-\widetilde{\vW})}\\
	 =& -\|T_{1}+T_{2}\|_{\overbar{\vW}^{-1}}^{2}-2\dotp{\overbar{\vW}^{-\frac{1}{2}}(T_{1}+T_{2}),\overbar{\vW}^{-\frac{1}{2}}T_{3}} -\|T_{3}\|_{\overbar{\vW}^{-1}}^{2}\\
		& \textstyle-\frac{4\theta-3}{4(1-\theta)}\|\vx-\vx^{+}\|^{2}_{(1-\theta)(\vI-\widetilde{\vW})}\\
\leq& -\|T_{1}+T_{2}\|_{\overbar{\vW}^{-1}}^{2}-2\dotp{\overbar{\vW}^{-\frac{1}{2}}(T_{1}+T_{2}),\overbar{\vW}^{-\frac{1}{2}}T_{3}}-(1+r_{1})\|T_{3}\|^{2}_{\overbar{\vW}^{-1}} \\
\leq& \textstyle-\frac{r_{1}}{1+r_{1}}\|T_{1}+T_{2}\|_{\overbar{\vW}^{-1}}^{2},
	\end{aligned}
\end{equation*}
where the last inequality comes from the Cauchy inequality $$\textstyle-2\dotp{a,b}\leq {1\over 1+r_1}\|a\|^{2}+(1+r_1)\|b\|^{2}.$$

Combining it with the last term $-\|\vx^{+}-\vx^{*}\|^{2}_{\vG}$ on the right hand side of~\eqref{thmeq3}, we have
\begin{align}
		&\textstyle -\frac{r_{1}}{1+r_{1}}\|T_{1}+T_{2}\|_{\overbar{\vW}^{-1}}^{2}-\|\vx^{+}-\vx^{*}\|^{2}_{\vG}\nonumber \\
\leq&\textstyle -\frac{r_{1}}{1+r_{1}}\|T_{2}\|_{\overbar{\vW}^{-1}}^{2}-\frac{2r_{1}}{1+r_{1}}\dotp{\overbar{\vW}^{-\frac{1}{2}}T_{1},\overbar{\vW}^{-\frac{1}{2}} T_{2}} -\frac{r_{1}}{1+r_{1}}\|T_{1}\|^2_{\overbar{\vW}^{-1}}\nonumber\\
&\textstyle-r_{2}\|T_{1}\|^{2}_{\overbar{\vW}^{-1}} -\frac{1}{2}\|\vx^{+}-\vx^{*}\|^{2}_{\vG}\nonumber\\
\leq&\textstyle -r_{3}\|\vy^{+}-\vy^{*}\|^{2}_{\overbar{\vW}^{-1}}-\frac{\xi}{2}\|\vx^{+}-\vx^{*}\|^{2}_{\vG}, \label{thmeq5}
	\end{align}
where $\xi<1$ is a small positive parameter, and $r_{2}$ and $r_{3}$ are defined as~\eqref{r2} and~\eqref{r3}, respectively. Since $\vG=(\vI-\vW)-2(\widetilde{\vW}-\vW),$ we have
\begin{equation}\label{decompose}
	\|\vx^{+}-\vx^{*}\|^{2}_{\vG}=\|\vx^{+}-\vx^{*}\|^{2}_{\vI-\vW}-2\|\vy-\vy^+\|^{2}_{\vM}.
\end{equation}
Therefore
\begin{align}
    &\textstyle -\frac{r_{1}}{1+r_{1}}\|T_{1}+T_{2}\|_{\overbar{\vW}^{-1}}^{2}-\|\vx^{+}-\vx^{*}\|^{2}_{\vG}\nonumber \\
\leq&\textstyle -r_{3}\|\vy^{+}-\vy^{*}\|^{2}_{\overbar{\vW}^{-1}}-\frac{\xi}{2}\|\vx^{+}-\vx^{*}\|^{2}_{\vI-\vW}  -\xi \|\vy-\vy^+\|_\vM^2\nonumber\\
\leq&\textstyle -r_{3}\|\vy^{+}-\vy^{*}\|^{2}_{\overbar{\vW}^{-1}}-\frac{\xi}{2}\|\vx^{+}-\vx^{*}\|^{2}_{\vI-\vW} +2\xi \|\vy^+-\vy^*\|_\vM^2+2\xi \|\vy-\vy^*\|_\vM^2\nonumber\\
\leq&\textstyle -(r_{3}/\lambda_{\max}(\overbar{\vW}\vM) -2\xi)\|\vy^{+}-\vy^{*}\|^{2}_{\vM}  -\frac{\xi}{2}\|\vx^{+}-\vx^{*}\|^{2}_{\vI-\vW}+2\xi \|\vy-\vy^*\|_\vM^2. \label{eq:ub2}
\end{align}
Let $\xi<  r_{3}/(4\lambda_{\max}(\overbar{\vW}\vM))$, then we have
$r_{3}/\lambda_{\max}(\overbar{\vW}\vM) -2\xi> 2\xi$.
Putting~\eqref{eq:ub1} and~\eqref{eq:ub2} together onto~\eqref{thmeq3}, we have
\begin{align*}
		&\textstyle \|\vx^{+}-\vx^{*}\|^{2}_{\vH}+\frac{\xi}{2}\|\vx^{+}-\vx^{*}\|^{2}_{\vI-\vW} +(1+(r_{3}/\lambda_{\max}(\overbar{\vW}\vM) -2\xi))\|\vy^{+}-\vy^{*}\|_{\vM}^{2}\\
\leq&\textstyle \big(1-\big(2\alpha-\frac{\alpha^{2}L}{\lambda_{\min}(\overbar{\vW})}\big){\mu_g}\big)\|\vx-\vx^{*}\|^{2}_{\vH} +\big(2\alpha-\frac{\alpha^{2}L}{\lambda_{\min}(\overbar{\vW})}\big)\eta\|\vx-\vx^*\|_{\vI-\vW}^2\\
	  & +(1+2\xi) \|\vy-\vy^*\|_\vM^2.
\end{align*}
Let $\rho$ be defined as~\eqref{rho}, we get~\eqref{linineq}.
Note that the choice of $\xi$ and $\eta$ affects the definition of $\vP$ and $\vQ$, but not the algorithm.
Hence for any $\alpha\in(0,\frac{2\lambda_{\min}(\overbar{\vW})}{L}),$ Q-linear convergence is guaranteed for $(\vx^{k}-\vx^{*},\vy^{k}-\vy^{*}). $

Because $\|\vx^{k}-\vx^{*}\|^2_{\vP}\leq\|\vx^{k}-\vx^{*}\|^2_{\vP}+\|\vy^{k}-\vy^{*}\|^{2}_{\vQ},$
the sequence $\{\|\vx^{k}-\vx^{*}\|^2_{\vP}\}_{k=1}^\infty$ converges R-linearly to $0$ at the rate of $\rho$. \qed
\end{proof}

Two special cases are not covered by the theorem: $\theta=1$ and $\widetilde{\vW}=\frac{\vI+\vW}{2}$.
When $\theta=1$, we have $r_{1}=\infty$ and $r_3={r_2\over1+r_2}$.
When $\widetilde{\vW}=\frac{\vI+\vW}{2}$, i.e., $\vG=\vzero$, we have $r_{2}=\infty$ and $r_{3}=\frac{r_{1}}{1+r_{1}}$. In both cases, the linear convergence rate is 
\begin{equation}\label{rho_sp}
    \begin{aligned}
    \rho=\max& \textstyle \Big\{\ 1-\Big(2\alpha-\frac{2\alpha^{2}L}{2-\theta+\theta\lambda_{\min}(\vW)}\Big)\mu_{\vg}, \Big(4\alpha-\frac{4\alpha^{2}L}{2-\theta+\theta\lambda_{\min}(\vW)}\Big)\frac{\eta}{\xi} ,\\
    & \textstyle 1-\frac{\beta r_{3}-4\xi(2-\theta\beta)}{\beta r_{3}+(1-2\xi)(2-\theta\beta)}\Big\},
    \end{aligned}
    \end{equation}
    where $\beta=1-\lambda_{2}(\vW)$ is the spectral gap. It is exactly the limit of $\rho$ in~\eqref{rho} with $r_1$ or $r_2$ approaching infinity.

\begin{remark}
The upper bound for the stepsize $\alpha$,
$\textstyle 
2(1-\theta+\theta\lambda_{\min}(\widetilde{\vW}))/L,$
is much larger than that in~\cite{Shi2015} for ensuring linear convergence, ${2\mu_{\vg}\lambda_{\min}(\widetilde{\vW})}/{L^{2}}$, when $\widetilde{\vW}$ is positive definite.
In the special case $\widetilde{\vW}=({\vI+\vW})/{2}$, we have $\alpha<(2-\theta+\theta\lambda_{\min}(\vW))/{L}$.
Since we can choose $\theta$ as close as possible to ${3}/{4}$, the upper bound of $\alpha$ attains $(3\lambda_{\min}(\vW)+5)/({4L})$, which coincides the optimal bound given in~\cite{li2017aprimal} for general convex functions.
In~\cite{li2017aprimal}, the linear convergence was shown under the strong convexity of all functions $\{f_i\}_{i=1}^n$.
\end{remark}

\subsection{NIDS without non-smooth term}
We consider NIDS next. In the smooth case, NIDS can be recovered by PAPC applied to the primal form of the decentralized consensus problem $$\Min_\vx~\vf(\vx),\quad\st \sqrt{\vI-\vW}\vx=\vzero.$$ It motivates us to show the inequality in Lemma~\ref{lemkeyineqnids}.

\cite[Lemma 1]{li2019decentralized} shows that, with the initialization $(\vd^0=\vzero,~\vx^0)$, the fixed point $(\vd^{*}\in\Range(\vI-\vW),\vx^{*})$ of~\eqref{reform_nids} satisfies 
\begin{subequations}\label{fixnid}
  \begin{align}
  \vd^{*}+\nabla\vf(\vx^{*})&= \vzero\label{fixnid:a},\\
  (\vI-\vW)\vx^{*}&= \vzero\label{fixnid:b},
  \end{align}
\end{subequations}
and $\vx^{*}$ is the consensual solution to the problem~\eqref{pbm}.
We will use the following important equality, which can be derived from~\eqref{reform_nids}
\beq\label{eq:dxconnection}
\begin{aligned}
\textstyle(\vI-{\vI-\vW\over 2})(\vd^{k+1}-\vd^k)=\frac{\vI-\vW}{2\alpha}(\vx^{k+1}-\vx^{*}).
\end{aligned}
\eeq

Motivated by the proof of EXTRA, we introduce another matrix to measure the distance to the fixed point.
We still pick $\theta\in(\frac{3}{4},1]$ such that
\beq\textstyle \theta\Big(\frac{\vI+\vW}{2}\Big)+(1-\theta)\vI=\vI-\theta\Big(\frac{\vI-\vW}{2}\Big)\succ\vzero.\eeq
Define a new symmetric matrix
\beq
\textstyle\widetilde{\vM}= 2(\vI-\vW)^{\dagger}-\theta\vI= \Big(\frac{\vI-\vW}{2}\Big)^{\dagger}-\theta\vI. \label{eq:Mtilde_NIDS}
\eeq
Then $\widetilde{\vM}$ is a norm over $\Range(\vI-\vW)$.
Note that $\widetilde{\vM}$ is invertible because $\widetilde{\vM}\vone=-\theta\vone$.
In the following proofs, we use the same simplification $\vx$ and $\vx^+$.

\begin{lemma}[Equality]\label{lemeqnids}
Let $\{(\vd^{k},\vx^{k})\}_{k=1}^{\infty}$ be the sequence generated by~\eqref{reform_nids}, we have the following two equalities:
\begin{subequations}\label{eqnid}
\begin{align}
 & \dotp{\vx^{k+1}-\vx^{*},\vd^{k+1}-\vd^{*}}= \alpha \dotp{\vd^{k+1}-\vd^{k},\vd^{k+1}-\vd^{*}}_{\widetilde{\vM}-(1-\theta)\vI} \label{eqnid:a}\\
 & \dotp{\vx^{k+1}-\vx^{*},\vd^{k+1}-\vd^{k}}= \alpha\|\vd^{k+1}-\vd^{k}\|^{2}_{\widetilde{\vM}-(1-\theta)\vI}\label{eqnid:b}.
\end{align}
\end{subequations}
\end{lemma}
\begin{proof}
Since $\vd^+-\vd^*\in\Range(\vI-\vW)$, we have
\begin{align}
  \dotp{\vx^{+}-\vx^{*},\vd^{+}-\vd^{*}}=& \textstyle \dotp{(\vI-\vW)(\vx^{+}-\vx^{*}),({\vI-\vW})^{\dagger}(\vd^{+}-\vd^{*})}\nonumber\\
=& \alpha\textstyle \dotp{(2\vI - (\vI-\vW))(\vd^{+}-\vd),({\vI-\vW})^{\dagger}(\vd^{+}-\vd^{*})}\nonumber\\
=& \alpha\textstyle \dotp{(2 (\vI-\vW)^{\dagger} - \vI )(\vd^{+}-\vd),\vd^{+}-\vd^{*}},\label{eq:eq_NIDS}
\end{align}
where the second equality follows~\eqref{eq:dxconnection}. 
Replacing $\vd^{*}$ with $\vd$ in~\eqref{eq:eq_NIDS}, we get~\eqref{eqnid:b} in the same way. \qed
\end{proof}

\begin{lemma}[A key inequality for NIDS]\label{lemkeyineqnids}
Let $\{(\vd^{k},\vx^{k})\}_{k=1}^{\infty}$ be the sequence generated by~\eqref{reform_nids}.
We have, with any $r_{4}\in(0,\theta-\frac{3}{4})$,
\begin{align}\label{ineqnids}
    & \|\vx^{k+1}-\vx^{*}\|^{2}+\alpha^2\|\vd^{k+1}-\vd^{*}\|^{2}_{\widetilde{\vM}+(\theta-\frac{1}{2}+2r_{4})\vI}\nonumber\\
\leq& \|\vx^{k}-\vx^{*}\|^{2}+\alpha^2\|\vd^{k}-\vd^{*}\|^{2}_{\widetilde{\vM}+(\theta-\frac{1}{2}-2r_{4})\vI} -\alpha^2\|\vd^{k}-\vd^{k+1}\|^{2}_{\widetilde{\vM}+(\theta-\frac{3}{4}-r_{4})\vI}\nonumber\\
    & +\alpha^2\|\nabla f(\vx^k)-\nabla f(\vx^*)\|^2  -2\alpha\dotp{\vx^{k}-\vx^{*},\nabla\vf(\vx^{k})-\nabla\vf(\vx^{*})}.
\end{align}
\end{lemma}

\begin{proof}
The iteration~\eqref{reform_nids} and the definition of $\widetilde{\vM}$ in~\eqref{eq:Mtilde_NIDS} show
\begin{equation*}
\begin{aligned}
 & 2\alpha\dotp{\vx-\vx^{*},\nabla\vf(\vx)-\nabla\vf(\vx^{*})}\\
=& 2\dotp{\vx-\vx^{*},\vx-\vx^{+}} -2\alpha\dotp{\vx-\vx^{*},\vd^{+}-\vd^{*}}\\
=& 2\dotp{\vx-\vx^{+},\vx-\vx^{*}} -2\alpha\dotp{\vx-\vx^+,\vd^{+}-\vd^{*}}-2\alpha\dotp{\vx^{+}-\vx^{*},\vd^{+}-\vd^{*}}\\
=& 2\dotp{\vx-\vx^{+},\vx-\alpha \vd^+-\vx^{*}+\alpha \vd^*} +2\alpha^2\dotp{\vd-\vd^+,\vd^{+}-\vd^{*}}_{\widetilde{\vM}-(1-\theta) \vI}\\
=& 2\dotp{\vx-\vx^{+},\vx^+-\vx^*+\alpha \nabla f(\vx)-\alpha \nabla f(\vx^*)} +2\alpha^2\dotp{\vd-\vd^+,\vd^{+}-\vd^{*}}_{\widetilde{\vM}-(1-\theta) \vI},
\end{aligned}
\end{equation*}
where the first and the last equalities use~\eqref{reform_nids:b} and the third one follows~\eqref{eqnid:a}. From~\eqref{reform_nids:b}, we obtain 
\begin{align}\label{themeqnids1}
 & 2\alpha\dotp{\vx-\vx^+,\nabla\vf(\vx)-\nabla\vf(\vx^{*})}\nonumber\\
=& \|\vx-\vx^+\|^{2}+\alpha^{2}\|\nabla\vf(\vx)-\nabla\vf(\vx^{*})\|^{2} -\|\vx-\vx^+-\alpha\nabla\vf(\vx)+\alpha\nabla\vf(\vx^{*})\|^{2}\nonumber\\
=& \|\vx-\vx^+\|^{2}+\alpha^{2}\|\nabla\vf(\vx)-\nabla\vf(\vx^{*})\|^{2}-\alpha^{2}\|\vd^{+}-\vd^{*}\|^{2}.
\end{align}
Together with the basic equality
$
2\dotp{a-b,b-c}  =\|a-c\|^{2}-\|b-c\|^{2}-\|a-b\|^{2},
$
we get
\begin{align}\label{eq3}
 & \|\vx^{+}-\vx^{*}\|^{2}+\alpha^2\|\vd^{+}-\vd^{*}\|^{2}_{\widetilde{\vM}-(1-\theta)\vI}\nonumber\\
=& \|\vx-\vx^{*}\|^{2}+\alpha^2\|\vd-\vd^{*}\|^{2}_{\widetilde{\vM}-(1-\theta)\vI}-\alpha^2\|\vd-\vd^{+}\|^{2}_{\widetilde{\vM}-(1-\theta)\vI}-\alpha^{2}\|\vd^{+}-\vd^{*}\|^{2}\nonumber\\
 & +\alpha^{2}\|\nabla\vf(\vx)-\nabla\vf(\vx^{*})\|^{2} -2\alpha\dotp{\vx-\vx^{*},\nabla\vf(\vx)-\nabla\vf(\vx^{*})}.
\end{align}

Since $r_{4}<\theta-\frac{3}{4}\leq 1/4$, the following inequality holds,
\begin{equation*}
\begin{aligned}
		\textstyle -(\frac{1}{2}-2r_{4})\|\vd^{+}-\vd^{*}\|^{2}
\leq\textstyle (\frac{1}{2}-2r_{4})\|\vd-\vd^{*}\|^{2}-(\frac{1}{4}-r_{4})\|\vd-\vd^{+}\|^{2}.
\end{aligned}
\end{equation*}
Adding it onto both sides of~\eqref{eq3}, we get~\eqref{ineqnids}.\qed
\end{proof}

\begin{theorem}[Q-linear convergence for NIDS]\label{thmnids}
Under Assumptions~\ref{asm1}-\ref{LipResSC}, we define
\begin{align}
r_5  =& \max\Big(2, {(\lambda_{\max}(\vI-\vW)-2)^2\over 2-({3\over4}+r_4)\lambda_{\max}(\vI-\vW)}\Big).\label{r5}
\end{align}
For any stepsize $\alpha\in(0,\frac{2}{L})$, we choose $\eta\in(0,{1\over \alpha(2-\alpha L)r_{5}})$ and
define
\begin{align}\label{rhonids}
\rho_{3}=\max \textstyle \left\{1-\alpha(2-\alpha L)\mu_{\vg},\alpha(2-\alpha L)\eta r_{5}\right.,
                     \textstyle \left.1-\frac{4r_{4}}{2\lambda_{\max}((\vI-\vW)^{+})-\frac{1}{2}+2r_{4}}\right\}<1,
\end{align}
Then we have
\begin{equation}\label{qconvnids}
\begin{aligned}
    & \|\vx^{k+1}-\vx^{*}\|^{2}_{\vI+{\vI-\vW\over r_5}}+\alpha^2\|\vd^{k+1}-\vd^{*}\|^{2}_{\vQ}\\
\leq& \rho(\|\vx^{k}-\vx^{*}\|^{2}_{\vI+{\vI-\vW\over r_5 }}+\alpha^2\|\vd^{k}-\vd^{*}\|^{2}_{\vQ}),
\end{aligned}
\end{equation}
where
$
\vQ\coloneqq  \textstyle\widetilde{\vM}+(\theta-\frac{1}{2}+2r_{4})\vI\succ\vzero.
$
\end{theorem}

\begin{proof}
Given any $\alpha\in(0,\frac{2}{L}),$ we have
\begin{equation*}\label{thmineqnids3}
\begin{aligned}
		& \alpha^{2}\|\nabla\vf(\vx)-\nabla\vf(\vx^{*})\|^{2} -2\alpha\dotp{\vx-\vx^{*},\nabla\vf(\vx)-\nabla\vf(\vx^{*})}\\
\leq&  -\alpha(2-\alpha L)\dotp{\vx-\vx^{*},\nabla\vf(\vx)-\nabla\vf(\vx^{*})}\\
	 =& -\alpha(2-\alpha L)\dotp{\vx-\vx^{*},\nabla\vf(\vx)-\nabla\vf(\vx^{*})} -\alpha(2-\alpha L)\eta\|\vx-\vx^{*}\|_{\vI-\vW}^{2}\\
		& +\alpha(2-\alpha L)\eta\|\vx-\vx^{*}\|_{\vI-\vW}^{2}\\
\leq& -\alpha(2-\alpha L)\mu_{\vg}\|\vx-\vx^{*}\|^{2} +\alpha(2-\alpha L)\eta\|\vx-\vx^{*}\|_{\vI-\vW}^{2},
\end{aligned}
\end{equation*}
where the first inequality is from~\eqref{Lip3} and the second one uses (restricted) strong convexity~\eqref{Stc2}.
Together with~\eqref{ineqnids}, we have
\begin{align}\label{thmineqnids1}
		& \|\vx^{+}-\vx^{*}\|^{2}+\alpha^2\|\vd^{+}-\vd^{*}\|^{2}_{\widetilde{\vM}+(\theta-\frac{1}{2}+2r_{4})\vI}\nonumber\\
\leq& \|\vx-\vx^{*}\|^{2}+\alpha^2\|\vd-\vd^{*}\|^{2}_{\widetilde{\vM}+(\theta-\frac{1}{2}-2r_{4})\vI} -\alpha^2\|\vd-\vd^{+}\|^{2}_{\widetilde{\vM}+(\theta-\frac{3}{4}-r_{4})\vI}\nonumber\\
		& -\alpha(2-\alpha L)\mu_{\vg}\|\vx-\vx^{*}\|^{2}+\alpha(2-\alpha L)\eta\|\vx-\vx^{*}\|_{\vI-\vW}^{2},
\end{align}
The equality~\eqref{eq:dxconnection} gives
\begin{align}
& \|\vx^+-\vx^*\|^2_{\vI-\vW} \nonumber\\
= & \|(\vI-\vW)(\vx^+-\vx^*)\|^2_{(\vI-\vW)^{\dagger}}
=  \alpha^2 \|(2\vI-(\vI-\vW))(\vd^{+}-\vd)\|^2_{(\vI-\vW)^{\dagger}}\nonumber\\
= & \alpha^2 \|\vd-\vd^+\|^2_{(2\vI-(\vI-\vW))(\vI-\vW)^{\dagger}(2\vI-(\vI-\vW))}
=  \alpha^2 \|\vd-\vd^+\|^2_{4(\vI-\vW)^{\dagger}-4\vI+(\vI-\vW)}\nonumber\\
\leq & \alpha^2 r_5 \|\vd-\vd^+\|^2_{\widetilde{\vM}+(\theta-{3\over4}-r_4)\vI}, \label{eq:dxbound}
\end{align}
where the second equality follows~\eqref{eq:dxconnection}, the fourth equality comes from $\vd-\vd^+\in\Range(\vI-\vW)$, and the inequality holds with the definition of $r_5$ in~\eqref{r5}. Combing~\eqref{thmineqnids1} and~\eqref{eq:dxbound}, we derive
\begin{align}
		& \textstyle\|\vx^{+}-\vx^{*}\|^{2}+{1\over r_5}\|\vx^+-\vx^*\|_{\vI-\vW}^2 +\alpha^2\|\vd^{+}-\vd^{*}\|^{2}_{\widetilde{\vM}+(\theta-\frac{1}{2}+2r_{4})\vI}\nonumber\\
\leq& (1-\alpha(2-\alpha L)\mu_{\vg})\|\vx-\vx^{*}\|^{2} +\alpha(2-\alpha L)\eta\|\vx-\vx^{*}\|_{\vI-\vW}^{2}\nonumber\\
		& +\alpha^2\|\vd-\vd^{*}\|^{2}_{\widetilde{\vM}+(\theta-\frac{1}{2}-2r_{4})\vI}.
\end{align}

Let $\rho_3$ be defined as~\eqref{rhonids}, and we show~\eqref{qconvnids}. Meanwhile, the Q-linear convergence of $(\vd^{k},\vx^{k})$ implies the R-linear convergence of $\vx^{k}$.\qed
\end{proof}

This theorem shows that NIDS is still linearly convergent over a relaxed $\vW$ and keeps the network-independent stepsize, which attains $\frac{2}{L}$ practically.

\section[Numerical Experiments]{Numerical Experiments\footnote{The datasets generated during and/or analysed during the current study are available from the corresponding author on reasonable request.}} \label{sec:num}
In this section, we compare the performance of EXTRA and NIDS over the relaxed mixing matrices in the following two scenarios:
\begin{itemize}
  \item Comparison of Decentralized Gradient Descent (DGD), EXTRA, and NIDS with different stepsizes for a doubly stochastic matrix $\vW$.
  \item Comparison of EXTRA and NIDS with different stepsizes for a relaxed matrix $\vW$.
\end{itemize}

We consider the following decentralized sensing problem.
Each agent $i\in\{1,\cdots,n\}$ has its own private measured data $M_{i}\in\mathbb{R}^{m_{i}\times p}$ and $y_{i}\in\mathbb{R}^{m_{i}}$ based on the unknown common variable $x\in\mathbb{R}^p$.
Suppose that $y_{i}=M_{i}x+e_{i}$ with independently identically distributed random noise $e_{i}\in\mathbb{R}^{m_{i}}$.
The goal is to estimate $x$ cooperatively over the network, and the problem is
\begin{equation*}
  \Min_{x}\ \bar{f}(x)=\frac{1}{n}\sum_{i=1}^{n}\frac{1}{2}\|M_{i}x-y_{i}\|_{2}^{2}.
\end{equation*}
The data $\{M_i\}_{i=1}^n$ and $x$ are generated from Gaussian distribution.
We normalize each $M_{i}$ such that $\|M_i^\top M_i\|=10$, i.e., $L=10$.
In both scenarios, we set $n=10$, $p=5$, $\vx^{0}=\vzero$, and $\widetilde{\vW}=\frac{\vI+\vW}{2}$ for EXTRA. 

For the first scenario, we construct the matrix $\vW$ based on the Metropolis constant edge weight matrix in~\cite[\S 2.4]{Shi2015}. 
In this case $\widetilde{\vW}$ is positive definite, and we can set $\theta=\frac{3}{4}$.
Then $\overbar{\vW}=\frac{5\vI+3\vW}{8}$.
We implement EXTRA with three different stepsizes: $\alpha_{1}=\frac{(1+\lambda_{\min}(\vW))\mu_{\bar{f}}}{100}$ (the stepsize for linear convergence in~\cite{Shi2015}), $\alpha_{2}=\frac{1+\lambda_{\min}(\vW)}{10}$ (the stepsize for convergence only in~\cite{Shi2015}), and $\alpha_{3}=\frac{5+3\lambda_{\min}(\vW)}{40}$ (our largest stepsize).
For NIDS, the stepsize is set to $\alpha_{4}=\frac{1}{5}$ although it is the upper bound of the stepsize which is not attainable in our proof theoretically.

\begin{figure}[!ht]
     \centering
     \includegraphics[width=.48\textwidth]{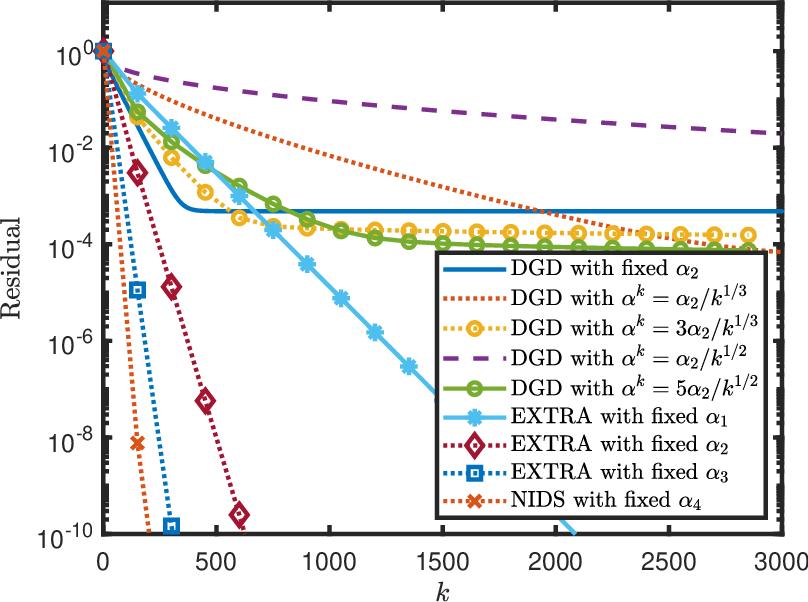}
     \includegraphics[width=.48\textwidth]{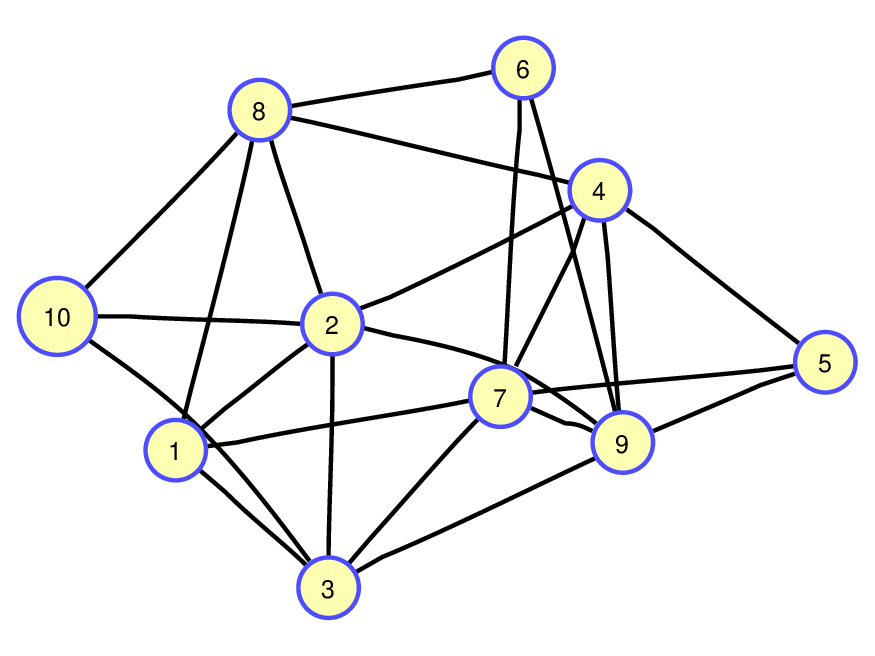}
     \caption{LEFT: the error $\frac{\|\vx^{k}-\vx^{*}\|_{F}}{\|\vx^{0}-\vx^{*}\|_{F}}$ vs iterations for DGD with different stepsizes, EXTRA with three stepsizes, and NIDS. RIGHT: The random network with 10 nodes.}\label{img1} 
\end{figure}

The result with $m_{i}=1$ is illustrated in Fig.~\ref{img1}.
Because we have $n>p$, the function $\bar f(x)$ is strongly convex with probability one. NIDS requires the least number of iteration to attain the expected tolerance.
Meanwhile, EXTRA with our proposed stepsize has better performance than that given in~\cite{Shi2015}.

\begin{figure}[!ht]
     \centering
     \includegraphics[width=.48\textwidth]{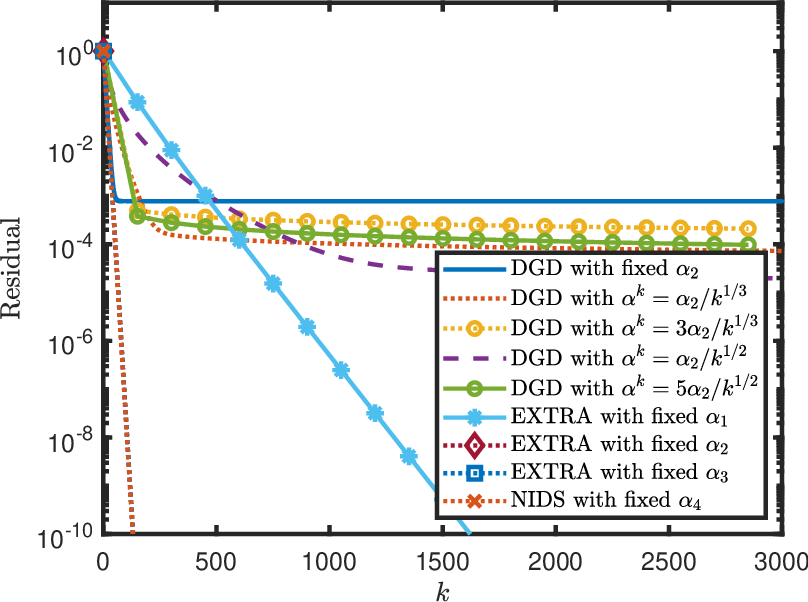}
     \includegraphics[width=.48\textwidth]{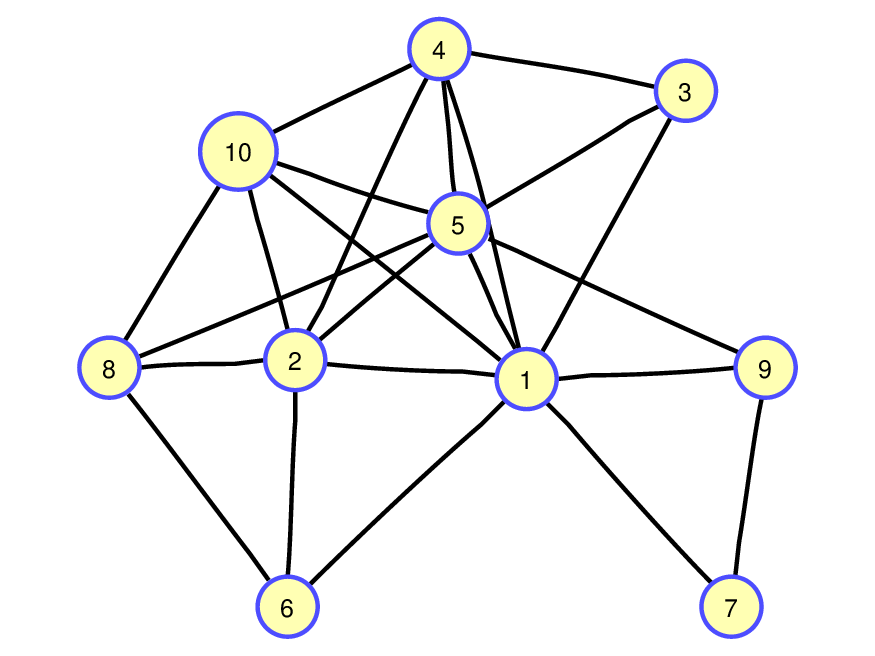}
     \caption{LEFT: the error $\frac{\|\vx^{k}-\vx^{*}\|_{F}}{\|\vx^{0}-\vx^{*}\|_{F}}$ vs iterations for DGD with different stepsizes, EXTRA with three stepsizes, and NIDS. RIGHT: The random network with 10 nodes.}\label{img2}
\end{figure}

\begin{figure}[!ht]
     \centering
     \includegraphics[width=.5\textwidth]{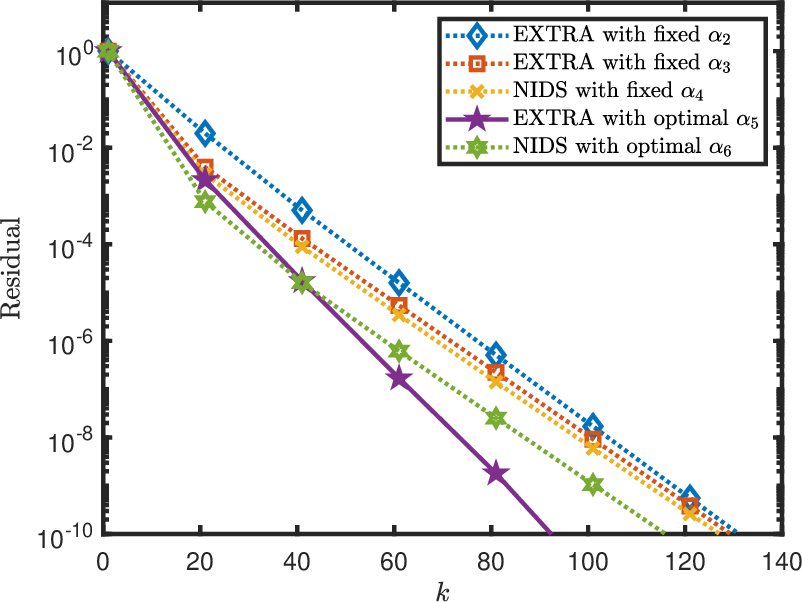}
     \caption{The comparison of proved stepsizes for EXTRA and NIDS with the optimal choice. 
     }\label{img4}
\end{figure}

Then, we set $m_{i}=10$ in Fig.~\ref{img2}.
In this case, individual functions $f_i(x)$ and $\vf(\vx)$ are strongly convex.
NIDS and EXTRA with the largest stepsize lead the performance. Here two results of EXTRA are the same as that of NIDS although they are set with different stepsizes. The observation may indicate that there is an optimal choice of stepsize between $\alpha_{2}$ and $\alpha_{3}$ for both EXTRA and NIDS. By setting $\alpha_{5}=\frac{5+3\lambda_{\min}(\vW)}{40+\mu_{\bar{f}}}$ for EXTRA and $\alpha_{6}=\frac{2}{10+\mu_{\bar{f}}}$, we compare these algorithms in Fig.~\ref{img4}. 
This figure suggests that the optimal stepsize may depends on the problem/functions. How to find the optimal stepsize is an important research topic and beyond the scope of this paper.

\begin{figure}[!ht]
     \begin{center}
     \includegraphics[width=.48\textwidth]{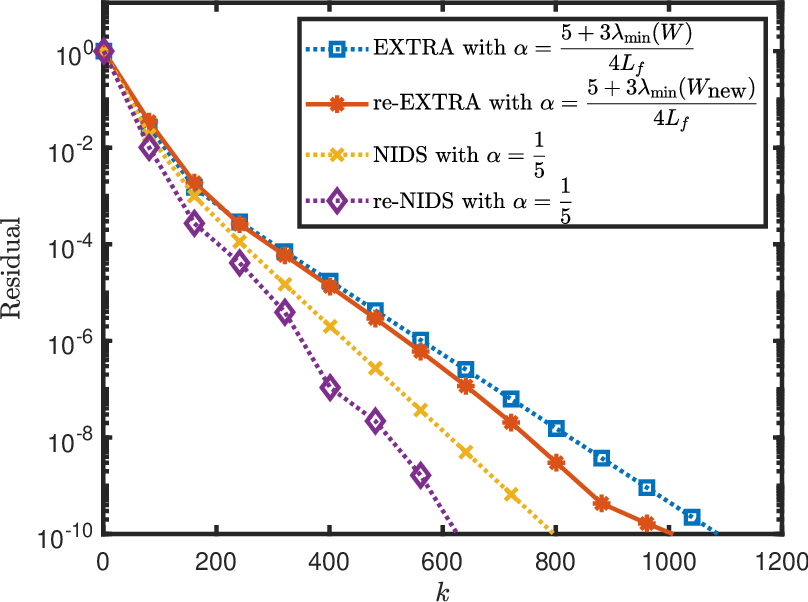}
     \includegraphics[width=.48\textwidth]{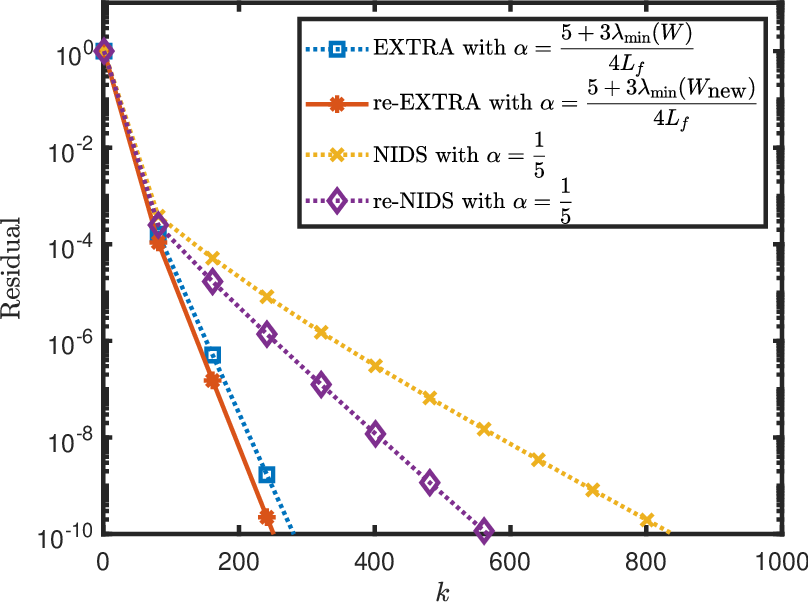}
     \end{center}
     \caption{The figure of residuals $\frac{\|\vx^{k}-\vx^{*}\|_{F}}{\|\vx^{0}-\vx^{*}\|_{F}}$ with respect to iteration. The left graph is for strongly convex $\bar{f}(x)$ and the right one is for strongly convex $\vf(\vx).$ re-EXTRA and re-NIDS stand for implementing EXTRA and NIDS over relaxed $\vW_{\mbox{new}}$.}\label{img3}
\end{figure}

Next, we turn to the relaxed mixing matrices.
Based on the previous created $\vW$, we replace it by $\vW_{\mbox{new}}={4\vW-\vI\over3}$ to scale the range of eigenvalues to $(-\frac{5}{3},1]$.
In this case, some diagonal entries of $\vW_{\mbox{new}}$ may be negative.
We consider the worst topology of network, line topology, i.e., each agent has at most two neighbors.
In this experiment, we solve the same problem using EXTRA and NIDS on the unrelaxed and relaxed mixing matrices, respectively, over the line. For NIDS, since the stepsize is network-independent, we relax the mixing matrix $\vW$ to $\vW_{\mbox{new}}$ more aggressively so that $\lambda_{\min}({\vW_{\mbox{new}}})$ approaches $-{5\over 3}$ and compare the performance with the unrelaxed case of NIDS under $\alpha=\frac{1}{5}.$ For EXTRA, we set the stepsize to $\alpha=\frac{5+3\lambda_{\min}(\vW)}{40}$, and compare the performance with the relaxed one under the stepsize $\alpha=\frac{5+3\lambda_{\min}(\vW_{\mbox{new}})}{40}$ where we only perturb $\vW$ mildly so that $\lambda_{\min}({\vW_{\mbox{new}}})$ approaches $-1$. The result is shown in Fig.~\ref{img3}.

From Fig.~\ref{img3}, if the topology of network is weak, switching to the relaxed mixing matrix may offer better performance when using NIDS and EXTRA to solve the problem. The improvement for NIDS is more distinguished.

\section{Conclusion}\label{sec:concl}
 In this paper, we relax the mixing matrices and prove the linear convergence of EXTRA and NIDS under the (restricted) strongly convexity assumption on $\bar{f}$. A larger upper bound of the stepsize is derived for EXTRA compared with that given in~\cite{Shi2015} and~\cite{alghunaim2019linearly}. 
 NIDS can choose a network-independent stepsize and this stepsize can be chosen as the same as that of centralized ones. We relax the conditions for the mixing matrices and the functions, while keeping the same stepsize. 

In numerical experiments on linear regression, EXTRA with the larger stepsize converges faster than using the $\mu_{\bar{f}}$-dependent stepsize in~\cite{Shi2015}. Over the unrelaxed mixing matrix, NIDS leads the performance in most cases and is the easiest to implement. If the topology of network is weak, using the relaxed mixing matrix can accelerate NIDS. For EXTRA, in general, we may not choose the mixing matrices to be relaxed due to the tiny improvement, but the larger stepsize derived in the relaxed case is competent to be considered.

\begin{acknowledgements}
This work is partially supported by NSF grants DMS-1621798 and DMS-2012439.
\end{acknowledgements}

\bibliographystyle{spmpsci_unsrt}
\bibliography{ref}

\end{document}